\nonstopmode
\documentclass[runningheads]{cl2emult}

\usepackage{amsmath}
\usepackage{amssymb}
\usepackage{bbm}          

\usepackage{url}          

\newcommand{\R}{\mathbbm{R}}
\newcommand{\Q}{\mathbbm{Q}}
\newcommand{\Z}{\mathbbm{Z}}
\newcommand{\N}{\mathbbm{N}}

\newcommand{\card}[1]{\lvert {#1} \rvert}
\newcommand{\setdef}[2]{\left\{ {#1} \, \vert \, {#2} \right\}}
\newcommand{\order}[1]{\mathcal{O}\!\left({#1}\right)}
\newcommand{\ordernorm}[1]{\mathcal{O}\!({#1})}
\newcommand{\graph}[1]{\mathcal{G}_{#1}}
\newcommand{\faclat}[1]{\mathcal{L}_{#1}}
\newcommand{\NP}{\ensuremath{\mathcal{NP}}}
\newcommand{\NumP}{\ensuremath{\#\mathcal{P}}}
\newcommand{\coNP}{\text{\rm co}\ensuremath{\mathcal{NP}}}
\DeclareMathOperator{\rank}{rank}
\DeclareMathOperator{\aff}{aff}
\DeclareMathOperator{\polylog}{polylog}
\DeclareMathOperator{\vol}{vol}

\newcounter{Problemcounter}

\newcommand{\mysection}[1]{\addcontentsline{ptc}{section}{#1}\section{#1}}

\makeatletter
\renewcommand{\subsection}{\@startsection{subsection}{2}{\z@}%
                          {-18\p@ \@plus -4\p@ \@minus -4\p@}%
                          {8\p@ \@plus 4\p@ \@minus 4\p@}%
                          {\normalfont\large\scshape
                           \rightskip=\z@ \@plus 8em\pretolerance=10000 }}
\makeatother

\makeatletter
\newcommand{\tableofproblems}{\@starttoc{ptc}}
\makeatother

\newenvironment{Problem}[1]%
{\refstepcounter{Problemcounter}
\begingroup
\subsection{\hspace{-3mm}#1}%
\endgroup
\addcontentsline{ptc}{subsection}{\arabic{Problemcounter}. \textsc{#1}}%
}%
{}

\newenvironment{ProblemIndex}[2]%
{\refstepcounter{Problemcounter}
\begingroup
\subsection{\hspace{-3mm}#1}%
\endgroup
\addcontentsline{ptc}{subsection}{\arabic{Problemcounter}. \textsc{#1}}%
}%
{}

\newenvironment{Definition}{\begin{list}{}%
{\setlength{\topsep}{0mm}%
\setlength{\partopsep}{0mm}%
\setlength{\parsep}{0mm}%
\setlength{\itemsep}{0mm}%
\setlength{\labelsep}{0mm}%
\settowidth{\labelwidth}{\textbf{Out}}%
\setlength{\leftmargin}{0mm}%
\addtolength{\leftmargin}{\labelwidth}%
\addtolength{\leftmargin}{\labelsep}%
\setlength{\itemindent}{0mm}%
}}{\end{list}}
\newcommand{\Input}{\item[\textbf{Input: }\hfill]}
\newcommand{\Output}{\item[\textbf{Output: }]}

\newcommand{\STATUS}[2]{\textbf{Status (general): } {#1}\\%
\textbf{Status (fixed dim.):} {#2}}

\newenvironment{Comment}{\begin{list}{}%
{\setlength{\topsep}{2mm}%
\setlength{\partopsep}{0mm}%
\setlength{\parsep}{0mm}%
\setlength{\itemsep}{0mm}%
\setlength{\labelsep}{0mm}%
\setlength{\labelwidth}{0mm}%
\setlength{\leftmargin}{0mm}%
\setlength{\itemindent}{0mm}%
}\item\relax}{\end{list}}
\newcommand{\Related}{\par\noindent\textbf{Related problems: }}

\title*{Some Algorithmic Problems\\in Polytope Theory}
\titlerunning{Some Algorithmic Problems in Polytope Theory}
\author{Volker Kaibel\thanks{Supported by the Deutsche
    Forschungsgemeinschaft, FOR~413/1--1 (Zi~475/3--1).} \and Marc E.
  Pfetsch} \authorrunning{Volker Kaibel and Marc E. Pfetsch}
\institute{TU~Berlin\\
  MA 6--2\\
  Stra\ss e des 17.~Juni~136\\
  10623~Berlin\\
  Germany\\
  \{\texttt{kaibel,pfetsch}\}\texttt{@math.tu-berlin.de} }

\begin{document}

\maketitle

\section{Introduction}

Convex polyhedra, i.e., the intersections of finitely many closed
affine half-spaces in $\R^d$, are important objects in various areas
of mathematics and other disciplines. In particular, the compact ones
among them (\emph{polytopes}), which equivalently can be defined as
the convex hulls of finitely many points in~$\R^d$, have been studied
since ancient times (e.g., the platonic solids). Polytopes appear as
building blocks of more complicated structures, e.g., in
(combinatorial) topology, numerical mathematics, or computer aided
design. Even in physics polytopes are relevant (e.g., in
crystallography or string theory).

Probably the most important reason for the tremendous growth of
interest in the theory of convex polyhedra in the second half of the
$20^{\text{th}}$ century was the fact that \emph{linear programming}
(i.e., optimizing a linear function over the solutions of a system of
linear inequalities) became a widespread tool to solve practical
problems in industry (and military). Dantzig's Simplex Algorithm,
developed in the late 40's, showed that geometric and combinatorial
know\-ledge of polyhedra (as the domains of linear programming problems)
is quite helpful for finding and analyzing solution procedures for
linear programming problems.

Since the interest in the theory of convex polyhedra to a large extent
comes from algorithmic problems, it is not surprising that many
algorithmic questions on polyhedra arose in the past. But also
inherently, convex polyhedra (in particular: polytopes) give rise to
algorithmic questions, because they can be treated as finite objects
by definition. This makes it possible to investigate (the smaller ones
among) them by computer programs (like the \texttt{polymake}-system
written by Gawrilow and Joswig, see~\cite{polyhome} and
\cite{GawJ00,GawJ01}). Once chosen to exploit this possibility, one
immediately finds oneself confronted with many algorithmic challenges.

This paper contains descriptions of 35~algorithmic problems about
polyhedra. The goal is to collect for each problem the current
knowledge about its computational complexity. Consequently, our
treatment is focused on theoretical rather than on practical subjects.
We would, however, like to mention that for many of the problems
computer codes are available.

Our choice of problems to be included is definitely influenced by
personal interest. We have not spent particular efforts to demonstrate
for each problem why we consider it to be relevant.  It may well be
that the reader finds other problems at least as interesting as the
ones we discuss. We would be very interested to learn about such
problems. The collection of problem descriptions presented in this
paper is intended to be maintained as a (hopefully growing)
list at \url{http://www.math.tu-berlin.de/~pfetsch/polycomplex/}.

Almost all of the problems are questions about polytopes. In some
cases the corresponding questions on general polyhedra are interesting
as well. It can be tested in polynomial time whether a polyhedron
specified by linear inequalities is bounded or not. This can be done
by applying Gaussian elimination and solving one linear program.

Roughly, the problems can be divided into two types: problems for
which the input are ``geometrical'' data and problems for which the
input is ``combinatorial'' (see below). Actually, it turned out that
it was rather convenient to group the problems we have selected into
the five categories ``Coordinate Descriptions''
(Sect.~\ref{sec:coordinate_descriptions}), ``Combinatorial Structure''
(Sect.~\ref{sec:combi}), ``Isomorphism''
(Sect.~\ref{sec:isomorphism}), ``Optimization''
(Sect.~\ref{sec:optimization}), and ``Realizability''
(Sect.~\ref{sec:realizability}).  Since the boundary complex of a
simplicial polytope is a simplicial complex, studying polytopes leads
to questions that are concerned with more general (polyhedral)
structures: \emph{simplicial complexes}.  Therefore, we have added a
category ``Beyond Polytopes'' (Sect.~\ref{sec:beyond}), where a few
problems concerned with general (abstract) simplicial complexes are
collected that are closely related to similar problems on polytopes.
We do not consider other related areas like oriented matroids.

The problem descriptions proceed along the following scheme.  First
input and output are specified. Then a summary of the knowledge on the
theoretical complexity is given, e.g., it is stated that the
complexity is unknown (``Open'') or that the problem is \NP-hard.
This is done for the case where the dimension (usually of the input
polytope) is part of the input as well as for the case of fixed
dimension; often the (knowledge on the) complexity status differs for
the two versions.  After that, comments on the problems are given
together with references. For each problem we tried to report on the
current state of knowledge according to the literature. Unless stated
otherwise, all results mentioned without citations are either
considered to be ``folklore'' or ``easy to prove.''  At the end
related problems in this paper are listed.

For all notions in the theory of polytopes that we use without
explanation we refer to Ziegler's book~\cite{Zie95}.  Similarly, for
the concepts from the theory of computational complexity that play a
role here we refer to Garey and Johnson's classical
text~\cite{GarJ79}. Whenever we talk about \emph{polynomial
  reductions} this refers to polynomial time Turing-reductions.  For
some of the problems the output can be exponentially large in the
input.  For these problems the interesting question is whether there
is a \emph{polynomial total time} algorithm, i.e., an algorithm whose
running time can be bounded by a polynomial in the sizes of the input
and the output (in contrast to a \emph{polynomial time} algorithm
whose running time would be bounded by a polynomial just in the input
size). Note that the notion of ``polynomial total time'' only makes
sense with respect to problems which explicitly require the output to
be non-redundant.

A very fundamental result in the theory of convex polyhedra is due to
Minkowski~\cite{Min96} and Weyl~\cite{Wey35}. For the special case of
polytopes (to which we restrict our attention from now on) it can be
formulated as follows. Every polytope $P\subset\R^d$ can be specified
by an $\mathcal{H}$- or by a $\mathcal{V}$-description. Here, an
$\mathcal{H}$-description consists of a finite set of linear
inequalities (defining closed affine half-spaces of~$\R^d$) such
that~$P$ is the set of all simultaneous solutions to these
inequalities. A $\mathcal{V}$-description consists of a finite set of
points in~$\R^d$ whose convex hull is~$P$.  If any of the two
descriptions is rational, then the other one can be chosen to be
rational as well. Furthermore, in this case the numbers in the second
description can be chosen such that their coding lengths depend only
polynomially on the coding lengths of the numbers in the first
description (see, e.g., Schrijver~\cite{Sch86}). In our context,
$\mathcal{H}$- and $\mathcal{V}$-descriptions are usually meant to be
rational. By linear programming, each type of description can be made
non-redundant in polynomial time (though it is unknown whether this is
possible in \emph{strongly} polynomial time, see
Problem~\ref{prob:geometric_lp}).

One of the basic properties of a polytope is its dimension. If the
polytope is given by a $\mathcal{V}$-description, then it can easily
be determined by Gaussian elimination (which, carefully done, is a
cubic algorithm; see, e.g., \cite{Sch86}). If the polyhedron is
specified by an $\mathcal{H}$-description, computing its dimension can
be done by linear programming (actually, this is polynomial time
equivalent to linear programming).

Furthermore, some of the problems may also be interesting in their
polar formulations, i.e., with ``the roles of $\mathcal{H}$- and
$\mathcal{V}$-descriptions exchanged.'' Switching to the polar
requires to have a relative interior point at hand, which is easy to
obtain if a $\mathcal{V}$-description is available, while it needs
linear programming if only an $\mathcal{H}$-description is specified.

We will especially be concerned with the combinatorial types of
polytopes, i.e., with their \emph{face lattices} (the sets of faces,
ordered by inclusion). In particular, some problems will deal with the
\emph{$k$-skeleton} of a polytope, which is the set of its faces of
dimensions less than or equal to~$k$, or with its \emph{$f$-vector},
i.e., the vector $(f_0(P),f_1(P),\dots,f_{d}(P))$, where $f_i(P)$ is
the number of $i$-dimensional faces (\emph{$i$-faces}) of the
$d$-dimensional polytope~$P$ (\emph{$d$-polytope}).  Talking of the
face lattice~$\faclat{P}$ of a polytope~$P$ will always refer to the
lattice as an abstract object, i.e., to \emph{any} lattice that is
isomorphic to the face lattice.  In particular, the lattice does not
contain any information on coordinates. Similarly, the
\emph{vertex-facet incidences} of~$P$ are given by any
matrix~$(a_{vf})$ with entries from~$\{0,1\}$, whose rows and columns
are indexed by the vertices and facets of~$P$, respectively, such
that~$a_{vf}=1$ if and only if vertex~$v$ is contained in facet~$f$.
Note that the vertex-facet incidences of a polytope completely
determine its face lattice.

A third important combinatorial structure associated with a
polytope~$P$ is its (abstract) graph~$\graph{P}$, i.e., any graph that
is isomorphic to the graph having the vertices of~$P$ as its nodes,
where two of them are adjacent if and only if their convex hull is a
(one-dimensional) face of~$P$. For simple polytopes, the (abstract)
graph determines the entire face lattice as well (see
Problem~\ref{prob:recsimp}). However, for general polytopes this is
not true.

Throughout the paper, $n$ refers to the number of vertices or points
in the given $\mathcal{V}$-description, respectively, depending on the
context.  Moreover, $m$ refers to the number of facets or inequalities
in the given $\mathcal{H}$-description, respectively, and $d$ refers
to the dimension of the polytope or the ambient space,
respectively.\smallskip

\noindent \textbf{Acknowledgment:} We thank the referee for many
valuable comments and G\"{u}nter M. Ziegler for carefully reading the
manuscript.


\mysection{Coordinate Descriptions}
\label{sec:coordinate_descriptions}

In this section problems are collected whose input are geometrical
data, i.e., the $\mathcal{H}$- or $\mathcal{V}$-description of a
polytope. Some problems which are also given by geometrical data
appear in Sections~\ref{sec:isomorphism} and \ref{sec:optimization}.

\begin{Problem}{Vertex Enumeration}
  \label{prob:vertex_enumeration}
  \begin{Definition}
    \Input Polytope $P$ in $\mathcal{H}$-description%
    \Output Non-redundant $\mathcal{V}$-description of $P$
  \end{Definition}
  \STATUS{Open; polynomial total time if $P$ is simple or
    simplicial}{Polynomial time}
  \begin{Comment}
    Let $d = \dim(P)$ and let $m$ be the number of inequalities in the
    input. It is well known that the number of vertices $n$ can be
    exponential ($\Omega (m^{\lfloor d/2 \rfloor})$) in the size of
    the input (e.g., Cartesian products of suitably chosen
    two-dimensional polytopes and prisms over them).
    
    \textsc{Vertex Enumeration} is strongly polynomially equivalent to
    Problem~\ref{prob:polytope_verification} (see Avis, Bremner, and
    Seidel~\cite{AviBS97}).  Since
    Problem~\ref{prob:facet_enumeration} is strongly polynomially
    equivalent to Problem~\ref{prob:polytope_verification} as well,
    \textsc{Vertex Enumeration} is also strongly polynomially
    equivalent to Problem~\ref{prob:facet_enumeration}.
    
    For fixed $d$, Chazelle~\cite{Cha93} found an $\order{m^{\lfloor
        d/2 \rfloor}}$ polynomial time algorithm, which is optimal by
    the Upper Bound Theorem of McMullen~\cite{McM70}. There exist
    algorithms which are faster than Chazelle's algorithm for small
    $n$, e.g., an $\order{m \log n + (m n)^{1-1/(\lfloor d/2 \rfloor
        +1)} \polylog m}$ algorithm of Chan~\cite{Cha96}.
    
    For general $d$, the reverse search method of Avis and Fukuda
    \cite{AviF92} solves the problem for \emph{simple} polyhedra in
    polynomial total time, using working space (without space for
    output) bounded polynomially in the input size. An algorithm of
    Bremner, Fukuda, and Marzetta \cite{BreFM98} solves the problem for
    \emph{simplicial} polytopes. Note that these algorithms need a
    vertex of $P$ to start from. Provan~\cite{Pro94} gives a
    polynomial total time algorithm for enumerating the vertices of
    polyhedra arising from networks.
    
    There are many more algorithms known for this problem -- none of
    them is a polynomial total time algorithm for general polytopes.
    See the overview article of Seidel~\cite{Sei97}. Most of these
    algorithms can be generalized to directly work for unbounded
    polyhedra, too. 
  \end{Comment}
  \Related \ref{prob:facet_enumeration},
  \ref{prob:polytope_verification},
  \ref{prob:face_lattice_enumeration}, \ref{prob:vertexnumber}
\end{Problem}

\begin{Problem}{Facet Enumeration}
  \label{prob:facet_enumeration}
  \begin{Definition}
    \Input Polytope $P$ in $\mathcal{V}$-description with $n$ points%
    \Output Non-redundant $\mathcal{H}$-description of $P$
  \end{Definition}
  \STATUS{Open; polynomial total time if $P$ is simple or
    simplicial}{Polynomial time}
  \begin{Comment}
    In~\cite{AviBS97} it is shown that \textsc{Facet Enumeration} is
    strongly polynomially equivalent to
    Problem~\ref{prob:polytope_verification} and thus to
    Problem~\ref{prob:vertex_enumeration} (see the comments there).
    
    For this problem, one can assume to have an interior point (e.g.,
    the vertex barycenter). \textsc{Facet Enumeration} is sometimes
    called the \emph{convex hull problem}.
  \end{Comment}
  \Related \ref{prob:vertex_enumeration},
  \ref{prob:polytope_verification}, \ref{prob:face_lattice_enumeration}
\end{Problem}

\begin{Problem}{Polytope Verification}
  \label{prob:polytope_verification}
  \begin{Definition}
    \Input Polytope $P$ given in $\mathcal{H}$-description, polytope
    $Q$ given in $\mathcal{V}$-description%
    \Output ``Yes'' if $P = Q$, ``No'' otherwise
  \end{Definition}
  \STATUS{Open; polynomial time if $P$ is simple or
    simplicial}{Polynomial time}
  \begin{Comment}
    \textsc{Polytope Verification} is strongly polynomially equivalent
    to Problem \ref{prob:vertex_enumeration} and Problem
    \ref{prob:facet_enumeration} (see the comments there).

    \textsc{Polytope Verification} is contained in $\coNP$: we can
    prove $Q\nsubseteq P$ by showing that some vertex of~$Q$ violates
    one of the inequalities describing~$P$. If $Q\subset P$ with
    $Q\neq P$ then there exists a point~$p$ of~$P\setminus Q$ with
    ``small'' coordinates (e.g., some vertex of~$P$ not contained
    in~$Q$) and a valid inequality for $Q$, which has ``small''
    coefficients and is violated by~$p$ (e.g., an inequality defining
    a facet of~$Q$ that separates~$p$ from~$Q$).  However, it is
    unknown whether \textsc{Polytope Verification} is in \NP.
    
    Since it is easy to check whether $Q \subseteq P$,
    \textsc{Polytope Verification} is
    Problem~\ref{prob:polytope_containment} restricted to the case
    that $Q \subseteq P$.
  \end{Comment}
  \Related \ref{prob:vertex_enumeration},
  \ref{prob:facet_enumeration}, \ref{prob:polytope_containment}
\end{Problem}

\begin{Problem}{Polytope Containment}
  \label{prob:polytope_containment}
  \begin{Definition}
    \Input Polytope $P$ given in $\mathcal{H}$-description, polytope
    $Q$ given in $\mathcal{V}$-description%
    \Output ``Yes'' if $P \subseteq Q$, ``No'' otherwise
  \end{Definition}
  \STATUS{$\coNP$-complete}{Polynomial time}
  \begin{Comment}
    Freund and Orlin~\cite{FreO85} proved that this problem is
    $\coNP$-complete. Note that the reverse question whether $Q
    \subseteq P$ is trivial. The questions where either both $P$ and
    $Q$ are given in $\mathcal{H}$-description or both in
    $\mathcal{V}$-description can be solved by linear programming
    (Problem \ref{prob:geometric_lp}), see Eaves and
    Freund~\cite{EavF82}. For fixed dimension, one can enumerate all
    vertices of $P$ in polynomial time (see Problem
    \ref{prob:vertex_enumeration}) and compare the descriptions of $P$
    and $Q$ (after removing redundant points).
  \end{Comment}
  \Related \ref{prob:polytope_verification}
\end{Problem}

\begin{Problem}{Face Lattice of Geometric Polytopes}
  \label{prob:face_lattice_enumeration}
  \begin{Definition}
    \Input Polytope $P$ in $\mathcal{H}$-description%
    \Output Hasse-diagram of the face lattice of $P$
  \end{Definition}
  \STATUS{Polynomial total time}{Polynomial time}
  \begin{Comment}
    See comments on Problem \ref{prob:vertex_enumeration}. Many
    algorithms for the \textsc{Vertex Enumeration Problem} in fact
    compute the whole face lattice of the polytope.
    Swart~\cite{Swa85}, analyzing an algorithm of Chand and
    Kapur~\cite{ChaK70}, proved that there exists a polynomial total
    time algorithm for this problem. For a faster algorithm see
    Seidel~\cite{Sei86}. Fukuda, Liebling, and Margot~\cite{FukLM97}
    gave an algorithm which uses working space (without space for the
    output) bounded polynomially in the input size, but it has to
    solve many linear programs. 
    
    For fixed dimension, the size of the output is polynomial in the
    size of the input; hence, a polynomial total time algorithm
    becomes a polynomial algorithm in this case.\smallskip
    
    The problem of enumerating the $k$-skeleton of $P$ seems to be
    open, even if $k$ is fixed. Note that, for fixed $k$, the latter
    problem can be solved by linear programming
    (Problem~\ref{prob:geometric_lp}) in polynomial time if the
    polytope is given in $\mathcal{V}$-description rather than in
    $\mathcal{H}$-description.
  \end{Comment}
  \Related \ref{prob:vertex_enumeration},
  \ref{prob:facet_enumeration}, \ref{prob:polytope_verification},
  \ref{prob:face_lattice}, \ref{prob:f_vector}
\end{Problem}

\begin{Problem}{Degeneracy Testing}
  \label{prob:degenerate}
  \begin{Definition}
    \Input Polytope $P$ in $\mathcal{H}$-description%
    \Output ``Yes'' if $P$ not simple (degenerate), ``No''
    otherwise
  \end{Definition}
  \STATUS{Strongly \NP-complete}{Polynomial time}
  \begin{Comment}
    Independently proved to be \NP-complete in the papers of
    Chandrasekaran, Kabadi, and Murty~\cite{ChaKM82} and Dyer
    \cite{Dye83}. Fukuda, Liebling, and Margot~\cite{FukLM97} proved
    that the problem is \emph{strongly} \NP-complete. For fixed
    dimension, one can enumerate all vertices in polynomial time (see
    Problem \ref{prob:vertex_enumeration}) and check whether they are
    simple or not.\smallskip
    
    Bremner, Fukuda, and Marzetta~\cite{BreFM98} noted that if $P$ is
    given in $\mathcal{V}$-description the problem is polynomial time
    solvable: enumerate the edges ($1$-skeleton, see Problem
    \ref{prob:face_lattice_enumeration}) and apply the Lower Bound
    Theorem.\smallskip
    
    Erickson~\cite{Eri99} showed that in the worst case
    $\Omega(m^{\lceil d/2 \rceil-1} + m \log m)$ sideness queries are
    required to test whether a polytope is simple. For odd $d$ this
    matches the upper bound. A \emph{sideness query} is a question of
    the following kind: given $d+1$ points $\vec{p}_0, \dots,
    \vec{p}_d$ in $\R^d$, does $\vec{p}_0$ lie ``above'', ``below'',
    or on the oriented hyperplane determined by $\vec{p}_1, \vec{p}_2,
    \dots, \vec{p}_d$.
  \end{Comment}
  \Related \ref{prob:vertex_enumeration}, \ref{prob:face_lattice_enumeration}
\end{Problem}

\begin{Problem}{Number of Vertices}
  \label{prob:vertexnumber}
  \begin{Definition}
    \Input Polytope $P$ in $\mathcal{H}$-description%
    \Output Number of vertices of $P$
  \end{Definition}
  \STATUS{\NumP-complete}{Polynomial time}
  \begin{Comment}
    Dyer~\cite{Dye83} and Linial~\cite{Lin86} independently proved
    that \textsc{Number of Vertices} is \NumP-complete. It follows
    that the problem of computing the $f$-vector of $P$ is \NumP-hard.
    Furthermore, Dyer~\cite{Dye83} proved that the decision version
    (``Given a number $k$, does $P$ have at least $k$ vertices?'') is
    strongly \NP-hard and remains \NP-hard when restricted to simple
    polytopes. It is unknown whether the decision problem is in \NP.
    
    If the dimension is fixed, one can enumerate all vertices in
    polynomial time (see Problem~\ref{prob:vertex_enumeration}).
  \end{Comment}
  \Related \ref{prob:vertex_enumeration}, \ref{prob:f_vector}
\end{Problem}

\begin{Problem}{Feasible Basis Extension}
  \label{prob:basisext}
  \begin{Definition}
    \Input Polytope $P$ given as $\{\vec{x} \in \R^s : A \vec{x} =
    \vec{b}, \vec{x} \geq \vec{0}\}$, a set $S \subseteq \{1, \dots, s\}$%
    \Output ``Yes'' if there is a feasible basis with an index set
    containing $S$, ``No'' otherwise
  \end{Definition}
  \STATUS{\NP-complete}{Polynomial time}
  \begin{Comment}
    See Murty~\cite{Mur72} (Garey and Johnson~\cite{GarJ79}, Problem
    [MP4]). For fixed dimension, one can enumerate all bases in
    polynomial time.
    
    The problem can be reformulated as follows. Let $P$ be defined by
    a finite set $\mathcal{H}$ of affine halfspaces and let $S$ be a
    subset of $\mathcal{H}$. Does $\bigcap \{H \in \mathcal{H} : H
    \notin S\}$ contain a vertex which is also a vertex of $P$?
  \end{Comment}
\end{Problem}

\begin{Problem}{Recognizing Integer Polyhedra}
  \begin{Definition}
    \Input Polytope $P$ in $\mathcal{H}$-description%
    \Output ``Yes'' if $P$ has only integer vertices, ``No''
    otherwise
  \end{Definition}
  \STATUS{Strongly $\coNP$-complete}{Polynomial time}
  \begin{Comment}
    The hardness-proof is by Papadimitriou and
    Yannakakis~\cite{PapY90}. For fixed dimension, one can enumerate
    all vertices (Problem~\ref{prob:vertex_enumeration}) and check
    whether they are integral in polynomial time.
  \end{Comment}
\end{Problem}

\begin{Problem}{Diameter}
  \begin{Definition}
    \Input Polytope $P$ in $\mathcal{H}$-description%
    \Output The diameter of $P$
  \end{Definition}
  \STATUS{$\NP$-hard}{Polynomial time}
  \begin{Comment}
    Frieze and Teng~\cite{FriT94} gave the proof of \NP-hardness. For
    fixed dimension, one can enumerate all vertices
    (Problem~\ref{prob:vertex_enumeration}), construct the graph and
    then compute the diameter in polynomial time. 
    
    The complexity status is unknown for simple polytopes. For
    simplicial polytopes the problem can be solved in polynomial time:
    Since simplicial polytopes have at most as many vertices as
    facets, one can enumerate their vertices (see
    Problem~\ref{prob:vertex_enumeration}), and finally compute the
    graph (and hence the diameter) from the vertex-facet incidences in
    polynomial time.
    
    If $P$ is given in $\mathcal{V}$-description, one can compute the
    graph ($1$-skeleton, see Problem
    \ref{prob:face_lattice_enumeration}) and hence the diameter in
    polynomial time.
  \end{Comment}
\end{Problem}

\begin{Problem}{Minimum Triangulation}
  \label{prob:triangulation}
  \begin{Definition}
    \Input Polytope $P$ in $\mathcal{V}$-description, positive integer $K$%
    \Output ``Yes'' if $P$ has a triangulation of size $K$ or less,
    ``No'' otherwise
  \end{Definition}
  \STATUS{\NP-complete}{\NP-complete}
  \begin{Comment}
    A \emph{triangulation} $\mathcal{T}$ of a $d$-polytope $P$ is a
    collection of $d$-simplices, whose union is $P$, their vertices
    are vertices of $P$, and any two simplices intersect in a common
    face (which might be empty). In particular, $\mathcal{T}$ is a
    (pure) $d$-dimensional geometric simplicial complex (see Section
    \ref{sec:beyond}). The size of $\mathcal{T}$ is the number of its
    $d$-simplices. Every (convex) polytope admits a triangulation.
    
    Below, De Loera, and Richter-Gebert~\cite{BelLR00b,BelLR00a}
    proved that \textsc{Minimum Triangulation} is \NP-complete for
    (fixed) $d\geq 3$. Furthermore, it is \NP-hard to compute a
    triangulation of minimal size for (fixed) $d\geq 3$.
  \end{Comment}
\end{Problem}

\begin{Problem}{Volume}
  \begin{Definition}
    \Input Polytope $P$ in $\mathcal{H}$-description%
    \Output The volume of $P$
  \end{Definition}
  \STATUS{\NumP-hard, FPRAS}{Polynomial time}
  \begin{Comment}
    Dyer and Frieze~\cite{DyeF88} showed that the general problem is
    \NumP-hard (and \NumP-easy as well). Dyer, Frieze, and
    Kannan~\cite{DyeFK91} found a \emph{Fully Polynomial Randomized
      Approximation Scheme} (\emph{FPRAS}) for the problem, i.e., a
    family $(A_{\varepsilon})_{\varepsilon>0}$ of randomized
    algorithms, where, for each $\varepsilon>0$, $A_{\varepsilon}$
    computes a number~$V_{\varepsilon}$ with the property that the
    probability of $(1-\varepsilon)\vol(P)\leq
    V_{\varepsilon}\leq(1+\varepsilon)\vol(P)$ is at
    least~$\frac{3}{4}$, and the running times of the algorithms
    $A_{\varepsilon}$ are bounded by a polynomial in the input size
    and~$\frac{1}{\varepsilon}$.
    \smallskip

    For fixed dimension, one can first compute all vertices of~$P$
    (see Problem~\ref{prob:vertex_enumeration}) and its face lattice
    (see Problem~\ref{prob:face_lattice_enumeration}) both in
    polynomial time.  Then one can construct some triangulation (see
    Problem~\ref{prob:triangulation}) of~$P$ (e.g., its barycentric
    subdivision) in polynomial time and compute the volume of~$P$ as
    the sum of the volumes of the (maximal) simplices in the
    triangulation.
    
    The complexity status of the analogue problem with the polytope
    specified by a $\mathcal{V}$-description is the same.
  \end{Comment}
\end{Problem}

\mysection{Combinatorial Structure}
\label{sec:combi}

In this section we collect problems that are concerned with computing
certain combinatorial information from compact descriptions of the
combinatorial structure of a polytope. Such compact encodings might be
the vertex-facet incidences, or, for simple polytopes, the abstract
graphs.  An example of such a problem is to compute the dimension of a
polytope from its vertex-facet incidences. Initialize a set~$S$ by the
vertex set of an arbitrary facet. For each facet~$F$ compute the
intersection of~$S$ with the vertex set of~$F$.  Replace~$S$ by a
maximal one among the proper intersections and continue. The dimension
is the number of ``rounds'' performed until~$S$ becomes empty.

All data is meant to be purely combinatorial. For all problems in this
section it is unknown if the ``integrity'' of the input data can be
checked, proved, or disproved in polynomial time. For instance, it is
rather unlikely that one can efficiently prove or disprove that a
lattice is the face lattice of some polytope (see
Problems~\ref{prob:steinitz_problem},
\ref{prob:simplicial_steinitz_problem}).

Sometimes, it might be worthwhile to exchange the roles of vertices
and facets by duality of polytopes.  Our choices of view points have
mainly been guided by personal taste.  \smallskip

Some orientations of the abstract graph~$\graph{P}$ of a simple
polytope~$P$ play important roles (although such orientations can also
be considered for non-simple polytopes, they have not yet proven to be
useful in the more general context). An orientation is called a
\emph{unique-sink orientation} (\emph{US-orient\-ation}) if it induces
a unique sink on every subgraph of~$\graph{P}$ corresponding to a
non-empty face of~$P$. A US-orientation is called an \emph{abstract
  objective function orientation} (\emph{AOF-orientation}) if it is
acyclic.  General US-orientations of graphs of cubes have recently
received some attention (Szab\'{o} and Welzl~\cite{SzaW01}).
AOF-orientations were used, e.g., by Kalai~\cite{Kal88}. Since their
linear extensions are precisely the shelling orders of the dual
polytope, they have been considered much earlier.

\begin{Problem}{Face Lattice of Combinatorial Polytopes}
  \label{prob:face_lattice}
  \begin{Definition}
    \Input Vertex-facet incidence matrix of a polytope~$P$
    \Output Hasse-diagram of the face lattice of~$P$
  \end{Definition}
  \STATUS{Polynomial total time}%
    {Polynomial time}
  \begin{Comment}
    Solvable in $\order{\min\{m,n\}\cdot\alpha\cdot\varphi}$ time,
    where~$m$ is the number of facets, $n$ is the number of vertices,
    $\alpha$ is the number of vertex-facet incidences, and~$\varphi$
    is the size of the face lattice~\cite{KaiP01}. Note that~$\varphi$
    is exponential in~$d$ (for fixed~$d$ it is polynomial in~$m$
    and~$n$).  Without (asymptotically) increasing the running time it
    is also possible to label each node in the Hasse diagram by the
    dimension and the vertex set of the corresponding face.

    It follows from~\cite{KaiP01} that one can compute the
    Hasse-diagram of the \emph{$k$-skeleton} (i.e., all faces of
    dimensions at most~$k$) of~$P$ in
    $\order{n\cdot\alpha\cdot\varphi^{\leq k}}$ time,
    where~$\varphi^{\leq k}$ is the number of faces of dimensions at
    most~$k$.  Since the latter number is in $\order{n^{k+1}}$, the
    $k$-skeleton can be computed in polynomial time (in the input
    size) for fixed~$k$.
  \end{Comment}
  \Related{\ref{prob:face_lattice_enumeration}, \ref{prob:f_vector}}
\end{Problem}

\begin{ProblemIndex}{$f$-Vector of Combinatorial Polytopes}{F-vector
    of combinatorial polytopes}
  \label{prob:f_vector}
  \begin{Definition}
    \Input Vertex-facet incidence matrix of a polytope~$P$
    \Output $f$-vector of~$P$
  \end{Definition}
  \STATUS{Open}%
  {Polynomial time}
  \begin{Comment}
    By the remarks on Problem~\ref{prob:face_lattice}, it is clear that 
    the first~$k$ entries of the $f$-vector can be computed in
    polynomial time for every fixed~$k$. 
    
    If the polytope is simplicial and a shelling (or a partition) of
    its boundary complex is available (see Problems~\ref{prob:AOF}
    and~\ref{prob:USO}), then one can compute the entire $f$-vector in
    polynomial time~\cite[Chap.~8]{Zie95}.
  \end{Comment}
  \Related{\ref{prob:vertexnumber}, \ref{prob:face_lattice}, 
  \ref{prob:AOF}, \ref{prob:USO}, \ref{prob:sc_f_vector}}
\end{ProblemIndex}

\begin{Problem}{Reconstruction of Simple Polytopes}
  \label{prob:recsimp}
  \begin{Definition}
    \Input The (abstract) graph~$\graph{P}$ of a simple polytope~$P$
    \Output The family of the subsets of nodes of~$\graph{P}$
    corresponding to the vertex sets of the facets of~$P$
  \end{Definition}
  \STATUS{Open}{Open}
  \begin{Comment}
    Blind and Mani~\cite{BM87} proved that the entire combinatorial
    structure of a simple polytope is determined by its graph. This is
    false for general polytopes (of dimension at least four), which is
    the main reason why we restrict our attention to simple polytopes
    for the remaining problems in this section.  Kalai~\cite{Kal88}
    gave a short, elegant, and constructive proof of Blind and Mani's
    result. However, the algorithm that can be derived from it has a
    worst-case running time that is exponential in the number of
    vertices of the polytope.  \smallskip
    
    In~\cite{JosKK01} it is shown that the problem can be formulated
    as a combinatorial optimization problem for which the problem to
    find an AOF-orientation of~$\graph{P}$ (see
    Problem~\ref{prob:AOF}) is strongly dual in the sense of
    combinatorial optimization. In particular, the vertex sets of the
    facets of~$P$ have a \emph{good characterization} in terms
    of~$\graph{P}$ (in the sense of Edmonds~\cite{Edm65}). The problem
    is polynomial time equivalent to computing the cycles
    in~$\graph{P}$ that correspond to the $2$-faces of~$P$.
    
    The problem can be solved in polynomial time in dimension at most
    three by computing a planar embedding of the graph, which can be
    done in linear time (Hopcroft and Tarjan~\cite{HopT74}, Mehlhorn and
    Mutzel~\cite{MehM96}).
  \end{Comment}
  \Related \ref{prob:facsimp}, \ref{prob:AOF}, \ref{prob:USO}
\end{Problem}

\begin{Problem}{Facet System Verification for Simple Polytopes}
  \label{prob:facsimp}
  \begin{Definition}
    \Input The (abstract) graph~$\graph{P}$ of a simple polytope~$P$
    and a family~$\mathcal{F}$ of subsets of nodes of~$\graph{P}$
    \Output ``Yes'' if~$\mathcal{F}$ is the family of subsets of nodes
    of~$\graph{P}$ that correspond to the vertex sets of the facets
    of~$P$, ``No'' otherwise
  \end{Definition}
  \STATUS{Open}{Open}
  \begin{Comment}
    In~\cite{JosKK01} it is shown that both the ``Yes''- as well as
    the ``No''-answer can be proved in polynomial time in the size
    of~$\graph{P}$ (provided that the integrity of the input data is
    guaranteed).  Any polynomial time algorithm for the construction
    of an AOF- or US-orientation (see Problems~\ref{prob:AOF}
    and~\ref{prob:USO}) of~$\graph{P}$ would yield a polynomial time
    algorithm for this problem (see~\cite{JosKK01}).
    
    Up to dimension three the problem can be solved in polynomial time
    (see the comments to Problems~\ref{prob:recsimp}
    and~\ref{prob:AOF}). 
  \end{Comment}
  \Related \ref{prob:recsimp}, \ref{prob:AOF}, \ref{prob:USO},
  \ref{prob:simplicial_steinitz_problem}
\end{Problem}

\begin{Problem}{AOF-Orientation}
  \label{prob:AOF}
  \begin{Definition}
    \Input The (abstract) graph~$\graph{P}$ of a simple polytope~$P$
    \Output An AOF-orientation of~$\graph{P}$
  \end{Definition}
  \STATUS{Open}{Open}
  \begin{Comment}
    (Simple) polytopes admit AOF-orientations because every linear
    function in general position induces an AOF-orientation.
    
    In~\cite{JosKK01} it is shown that one can formulate the problem
    as a combinatorial optimization problem, for which a strongly dual
    problem in the sense of combinatorial optimization exists (see the
    comments to Problem~\ref{prob:recsimp}).  Thus, the
    AOF-orientations of~$\graph{P}$ have a good characterization (see
    Problem~\ref{prob:facsimp}) in terms of~$\graph{P}$, i.e., there
    are polynomial size proofs for both cases an orientation being an
    AOF-orientation or not (provided that the integrity of the input
    data is guaranteed).  However, it is unknown if it is possible to
    check in polynomial time if a given orientation is an AOF-orientation.
    
    In dimensions one and two the problem is trivial. For a
    three-dimensional polytope~$P$ the problem can be solved in
    polynomial time, e.g., by producing a plane drawing of~$\graph{P}$
    with convex faces (see Tutte~\cite{Tut63}) and sorting the nodes
    with respect to a linear function (in general position).

    A polynomial time algorithm would lead to a polynomial algorithm for
    Problem~\ref{prob:facsimp} (see~\cite{JosKK01}).
    
    By duality of polytopes, the problem is equivalent to the problem
    of finding a shelling order of the facets of a simplicial polytope
    from the upper two layers of its face lattice. It is unknown
    whether it is possible to find in polynomial time a shelling order
    of the facets, even if the polytope is given by its entire face
    lattice. With this larger input, however, it is possible to check
    in polynomial time whether a given ordering of the facets is a
    shelling order.
  \end{Comment}
  \Related{\ref{prob:facsimp}, \ref{prob:USO}, \ref{prob:shellability}}
\end{Problem}

\begin{Problem}{US-Orientation}
  \label{prob:USO}
  \begin{Definition}
    \Input The (abstract) graph~$\graph{P}$ of a simple polytope~$P$
    \Output A US-orientation of~$\graph{P}$
  \end{Definition}
  \STATUS{Open}{Open}
  \begin{Comment}
    Since AOF-orientations are US-orientations, it follows from the
    remarks on Problem~\ref{prob:AOF} that (simple) polytopes admit
    US-orientations and that the problem can be solved in polynomial
    time up to dimension three.  By slight adaptions of the arguments
    given in~\cite{JosKK01}, one can prove that a polynomial time
    algorithm for this problem would yield a polynomial time algorithm
    for Problem~\ref{prob:facsimp} as well.
    
    In contrast to Problem~\ref{prob:AOF}, no good characterization of
    US-orien\-ta\-tions is known.
    
    It is not hard to see that, by duality of polytopes, the problem
    is equivalent to the problem of finding from the upper two layers
    a partition of the face lattice of a simplicial polytope into
    intervals whose upper bounds are the facets (i.e., a partition in
    the sense of Stanley~\cite{Sta79}). Similar to the situation with
    shelling orders, it is even unknown whether such a partition can
    be found in polynomial time if the polytope is specified by its
    entire face lattice. Again, with the entire face lattice as input
    it can be checked in polynomial time whether a family of subsets
    of the face lattice is a partition in that sense.
  \end{Comment}
  \Related{\ref{prob:facsimp}, \ref{prob:AOF}, \ref{prob:partitionability}}
\end{Problem}

\mysection{Isomorphism}
\label{sec:isomorphism}

Two polytopes $P_1\subset\R^{d_1}$ and $P_2\subset\R^{d_2}$ are
\emph{affinely equivalent} if there is a one-to-one affine map
$T:\aff(P_1)\longrightarrow\aff(P_2)$ between the affine hulls
of~$P_1$ and~$P_2$ with $T(P_1)=P_2$. Two polytopes are
\emph{combinatorially equivalent} (or \emph{isomorphic}) if their face
lattices are isomorphic. It is not hard to see that affine equivalence
implies combinatorial equivalence.

As soon as one starts to investigate structural properties of
polytopes by means of computer programs, algorithms for deciding
whether two polytopes are isomorphic become relevant.

Some problems in this section are known to be hard in the sense that
the \emph{graph isomorphism problem} can polynomially be reduced to
them. Although this problem is not known (and even not expected) to be
$\NP$-complete, all attempts to find a polynomial time algorithm for
it have failed so far. Actually, the same holds for a lot of problems
that can be polynomially reduced to the graph isomorphism problem
(see, e.g., Babai~\cite{Bab95}).

\begin{ProblemIndex}{Affine Equivalence of
    $\protect\mathcal{V}$-Polytopes}{Affine Equivalence of V-Polytopes}
  \label{prob:aff_equi}
  \begin{Definition}
    \Input Two polytopes~$P$ and~$Q$ given in
    $\mathcal{V}$-description%
    \Output ``Yes'' if $P$ is affinely equivalent to $Q$, ``No''
    otherwise
  \end{Definition}
  \STATUS{Graph isomorphism hard}{Polynomial time}
  \begin{Comment}
    The graph isomorphism problem can  polynomially be reduced to the
    problem of checking the affine equivalence of two
    polytopes~\cite{KaiS01}. The problem remains graph isomorphism
    hard if $\mathcal{H}$-descriptions are additionally provided as
    input data and/or if one restricts the input to simple or simplicial
    polytopes.
    
    For polytopes of bounded dimension the problem can be solved in
    polynomial time by mere enumeration of affine bases among the
    vertex sets.
  \end{Comment}
  \Related{\ref{prob:comb_equi}}
\end{ProblemIndex}

\begin{ProblemIndex}{Combinatorial Equivalence of $\protect\mathcal{V}$-Polytopes}{Combinatorial Equivalence of V-Polytopes}
  \label{prob:comb_equi}
  \begin{Definition}
    \Input Two polytopes~$P$ and~$Q$ given in
    $\mathcal{V}$-description%
    \Output ``Yes'' if $P$ is combinatorially equivalent to $Q$,
    ``No'' otherwise
  \end{Definition}
  \STATUS{\coNP-hard}{Polynomial time}
  \begin{Comment}
    Swart~\cite{Swa85} describes a reduction of the subset-sum problem
    to the negation of the problem.
    
    For polytopes of bounded dimension the problem can be solved in
    polynomial time (see Problems~\ref{prob:facet_enumeration} and~\ref{prob:isovfi}).
  \end{Comment}
  \Related{\ref{prob:facet_enumeration}, \ref{prob:aff_equi}, \ref{prob:isovfi}}
\end{ProblemIndex}

\begin{Problem}{Polytope Isomorphism}
  \label{prob:polyiso}
  \begin{Definition}
    \Input The face lattices~$\faclat{P}$ and~$\faclat{Q}$ of two
    polytopes~$P$ and~$Q$, respectively%
    \Output ``Yes'' if~$\faclat{P}$ is isomorphic to~$\faclat{Q}$,
    ``No'' otherwise
  \end{Definition}
  \STATUS{Open}{Polynomial time}
  \begin{Comment}    
    The problem can be solved in polynomial time in constant
    dimension (see Problem~\ref{prob:isovfi}). In general, the problem 
    can easily be reduced to the graph isomorphism problem
  \end{Comment}
  \Related{\ref{prob:isovfi}, \ref{prob:selfdual}}
\end{Problem}

\begin{Problem}{Isomorphism of vertex-facet incidences}
  \label{prob:isovfi}
  \begin{Definition}
    \Input Vertex-facet incidence matrices~$A_P$ and~$A_Q$ of 
    polytopes~$P$ and~$Q$, respectively%
    \Output ``Yes'' if~$A_P$ can be transformed into~$A_Q$ by row and
    column permutations, ``No'' otherwise
  \end{Definition}
  \STATUS{Graph isomorphism complete}{Polynomial time}
  \begin{Comment}
    The problem remains graph isomorphism complete even if
    $\mathcal{V}$- and~$\mathcal{H}$-descript\-ions of~$P$ and~$Q$ are
    part of the input data~\cite{KaiS01}.
    
    In constant dimension the problem can be solved in polynomial time
    by a reduction~\cite{KaiS01} to the graph isomorphism problem for
    graphs of bounded degree, for which a polynomial time algorithm is
    known (Luks~\cite{Luk82}).
    
    Problem~\ref{prob:polyiso} can  polynomially be reduced to this
    problem. For polytopes of bounded dimension both problems are
    polynomial time equivalent.
  \end{Comment}
  \Related{\ref{prob:polyiso}, \ref{prob:comb_equi}}
\end{Problem}

\begin{Problem}{Selfduality of Polytopes}
  \label{prob:selfdual}
  \begin{Definition}
    \Input Face Lattice of a polytope~$P$ 
    \Output ``Yes'' if~$P$ is isomorphic to its dual,
    ``No'' otherwise
  \end{Definition}
  \STATUS{Open}{Polynomial time}
  \begin{Comment}
    This is a special case of problem~\ref{prob:polyiso}. In
    particular, it is solvable in polynomial time in bounded
    dimensions. 
    
    It is easy to see that deciding whether a general 0/1-matrix~$A$
    (not necessarily a vertex-facet incidence matrix of a polytope)
    can be transformed into~$A^T$ by permuting its rows and columns is
    graph isomorphism complete.
  \end{Comment}
  \Related{\ref{prob:polyiso}}
\end{Problem}

\mysection{Optimization}
\label{sec:optimization}

In this section, next to the original linear programming problem, we
describe some of its relatives. In particular, combinatorial
abstractions of the problem are important with respect to polytope
theory (and, more general, discrete geometry). We pick out the aspect
of combinatorial cube programming here (and leave aside abstractions
like general combinatorial linear programming, LP-type problems, and
oriented matroid programming), since it has received considerable
attention lately.

\begin{Problem}{Geometric Linear Programming}
  \label{prob:geometric_lp}
  \begin{Definition}
    \Input $\mathcal{H}$-description of a polyhedron~$P\subset\Q^d$,
    $\vec{c} \in\Q^d$%
    \Output$\inf\setdef{\vec{c}^T \vec{x}}{\vec{x} \in
      P}\in\Q\cup\{-\infty,\infty\}$ and, if the infimum is finite, a
    point where the infimum is attained.
  \end{Definition}
  \STATUS{Polynomial time; no \emph{strongly} polynomial time
    algorithm known}{Linear time in~$m$ (the number of inequalities)}
  \begin{Comment}
    The first polynomial time algorithm was a variant of the
    \emph{ellipsoid algorithm} due to Khachiyan~\cite{Kha79}. Later,
    also \emph{interior point methods} solving the problem in
    polynomial time were discovered (Karmarkar~\cite{Kar84}).
    
    Megiddo found an algorithm solving the problem for a fixed
    number~$d$ of variables in $\order{m}$ arithmetic operations
    (Megiddo~\cite{Meg84}).
    
    No strongly polynomial time algorithm (performing a number of
    arithmetic operations that is bounded polynomially in~$d$ and the
    number of half-spaces rather than in the coding lengths of the
    input coordinates) is known.  In particular, no polynomial time
    variant of the \emph{simplex algorithm} is known. However, a
    randomized version of the simplex algorithm solves the problem in
    (expected) subexponential time (Kalai~\cite{Kal92}, Matou\v{s}ek,
    Sharir, and Welzl~\cite{MatSW96}).
  \end{Comment}
  \Related \ref{prob:optimalvertex}, \ref{prob:specifiedvalue},
  \ref{prob:AOFprog}
\end{Problem}

\begin{Problem}{Optimal Vertex}
  \label{prob:optimalvertex}
  \begin{Definition}
    \Input $\mathcal{H}$-description of a polyhedron~$P\subset\Q^d$,
    $\vec{c} \in \Q^d$%
    \Output$\inf\setdef{\vec{c}^T\vec{v}}{\vec{v}\text{ vertex of
        }P}\in\Q\cup\{\infty\}$ and, if the infimum is finite, a
    vertex where the infimum is attained.
  \end{Definition}
  \STATUS{Strongly \NP-hard}{Polynomial time}
  \begin{Comment}
    Proved to be strongly \NP-hard by Fukuda, Liebling, and
    Margot~\cite{FukLM97}. By linear programming this problem can be
    solved in polynomial time if~$P$ is a polytope. In fixed dimension
    one can compute all vertices of~$P$ in polynomial time (see
    Problem~\ref{prob:vertex_enumeration}).
  \end{Comment}
  \Related{\ref{prob:vertex_enumeration}, \ref{prob:geometric_lp},
    \ref{prob:specifiedvalue}}
\end{Problem}

\begin{Problem}{Vertex with specified objective value}
  \label{prob:specifiedvalue}
  \begin{Definition}
    \Input $\mathcal{H}$-description of a polyhedron~$P\subset\Q^d$,
    $\vec{c} \in \Q^d$, $C\in\Q$%
    \Output ``Yes'' if there is a vertex~$\vec{v}$ of $P$ with
    $\vec{c}^T \vec{v} = C$; ``No'' otherwise
  \end{Definition}
  \STATUS{Strongly \NP-complete}{Polynomial time}
  \begin{Comment}
    Proved to be \NP-complete by Chandrasekaran, Kabadi, and
    Murty~\cite{ChaKM82} and strongly \NP-complete by Fukuda,
    Liebling, and Margot~\cite{FukLM97}. The problem remains strongly
    \NP-complete even if the input is restricted to
    polytopes~\cite{FukLM97}.
  \end{Comment}
  \Related{\ref{prob:geometric_lp}, \ref{prob:optimalvertex}}
\end{Problem}

\begin{Problem}{AOF Cube Programming}
  \label{prob:AOFprog}
  \begin{Definition}
    \Input An oracle for a function
    $\sigma:\{0,1\}^d\longrightarrow\{+,-\}^d$ defining an
    AOF-orientation of the graph of the $d$-cube
    \Output The sink of the orientation
  \end{Definition}
  \STATUS{Open}{Constant time}
  \begin{Comment}
    The problem can be solved in a subexponential number of oracle
    calls by the \emph{random facet} variant of the simplex algorithm
    due to Kalai~\cite{Kal92}. For a derivation of the explicit
    bound~$e^{2\sqrt{d}}-1$ see G\"artner~\cite{Gae01}.

    In fixed dimension the problem is trivial  by mere enumeration.
    
    The problem generalizes linear programming problems whose sets of
    feasible solutions are combinatorially equivalent to cubes.
  \end{Comment}
  \Related{\ref{prob:geometric_lp}, \ref{prob:USprog}}
\end{Problem}

\begin{Problem}{USO Cube Programming}
  \label{prob:USprog}
  \begin{Definition}
    \Input An oracle for a function
    $\sigma:\{0,1\}^d\longrightarrow\{+,-\}^d$ defining a
    US-orientation of the graph of the $d$-cube
    \Output The sink of the orientation
  \end{Definition}
  \STATUS{Open}{Constant time}
  \begin{Comment}
    Szab\'{o} and Welzl~\cite{SzaW01} describe a randomized algorithm
    solving the problem in an expected number of
    $\ordernorm{\alpha^d}$ oracle calls with
    $\alpha=\sqrt{43/20}<1.467$ and a deterministic algorithm that
    needs $\order{1.61^d}$ oracle calls. Plugging an optimal algorithm
    for the three-dimensional case (found by G\"unter Rote) into their
    algorithm, Szab\'{o} and Welzl even obtain an~$\order{1.438^d}$
    randomized algorithm.

    The problem not only generalizes Problem~\ref{prob:AOFprog},
    but also certain linear complementary problems and smallest
    enclosing ball problems.

    In fixed dimension the problem is trivial by mere enumeration.
  \end{Comment}
  \Related{\ref{prob:AOFprog}}
\end{Problem}

\mysection{Realizability}
\label{sec:realizability}

In this section problems are discussed which bridge the gap from
combinatorial descriptions of polytopes to geometrical descriptions,
i.e., it deals with questions of the following kind: given
combinatorial data, does there exist a polytope which ``realizes''
this data? E.g., given a $0/1$-matrix is this the matrix of
vertex-facet incidences of a polytope? The problems of computing
combinatorial from geometrical data is discussed in
Section~\ref{sec:coordinate_descriptions}.

The problems listed in this section are among the first ones asked in
(modern) polytope theory, going back to the work of Steinitz and
Radermacher in the 1930's~\cite{SteR34}.

\begin{Problem}{Steinitz Problem}
  \label{prob:steinitz_problem}
  \begin{Definition}
    \Input Lattice $\mathcal{L}$%
    \Output ``Yes'' if $\mathcal{L}$ is isomorphic to the face lattice of a
    polytope, ``No'' otherwise
  \end{Definition}
  \STATUS{\NP-hard}{\NP-hard}
  \begin{Comment}
    If $\mathcal{L}$ is isomorphic to the face lattice of a polytope,
    it is ranked, atomic, and coatomic. These properties can be tested
    in polynomial time in the size of $\mathcal{L}$. Furthermore, in
    this case, the dimension $d$ of a candidate polytope has to be
    $\rank{\mathcal{L}}-1$.
    
    The problem is trivial for dimension $d \leq 2$.  Steinitz's
    Theorem allows to solve $d=3$ in polynomial time: construct the
    (abstract) graph $G$, test if the facets can consistently be
    embedded in the plane (linear time~\cite{HopT74,MehM96}) and check for
    $3$-connectedness (in linear time, see Hopcroft and
    Tarjan~\cite{HopT73}).
    
    Mn{\"e}v proved that the Steinitz Problem for $d$-polytopes with
    $d+4$ vertices is \NP-hard~\cite{Mne88}. Even more,
    Richter-Gebert~\cite{RiG96} proved that for (fixed) $d \geq 4$ the
    problem is \NP-hard.
    
    For fixed $d\geq 4$ it is neither known whether the problem is
    in~\NP\ nor whether it is in~$\coNP$. It seems unlikely to be
    in~\NP, since there are $4$-polytopes which cannot be realized
    with rational coordinates of coding length which is bounded by a
    polynomial in $\card{\mathcal{L}}$ (see
    Richter-Gebert~\cite{RiG96}).
  \end{Comment}
  \Related \ref{prob:simplicial_steinitz_problem}
\end{Problem}

\begin{Problem}{Simplicial Steinitz Problem}
  \label{prob:simplicial_steinitz_problem}
  \begin{Definition}
    \Input Lattice $\mathcal{L}$%
    \Output ``Yes'' if $\mathcal{L}$ is isomorphic to the face lattice
    of a simplicial polytope, ``No'' otherwise
  \end{Definition}
  \STATUS{\NP-hard}{Open}
  \begin{Comment}
    As for Problem \ref{prob:steinitz_problem}, $\mathcal{L}$ is
    ranked, atomic, and coatomic if the answer is ``Yes.'' In this
    case, the dimension $d$ of any matched polytope is
    $\rank{\mathcal{L}}-1$.
    
    As for general polytopes (Problem \ref{prob:steinitz_problem}),
    this problem is polynomial time solvable in dimension $d \leq 3$.
    
    The problem is \NP-hard, which follows from the above mentioned
    fact that the Steinitz problem for $d$-polytopes with $d+4$
    vertices is \NP-hard and a construction (Bokowski and
    Sturmfels~\cite{BokS89}) which generalizes it to the
    \emph{simplicial} case (but increases the dimension). It is,
    however, open whether the problem is \NP-hard for \emph{fixed}
    dimension. For fixed $d\geq 4$, it is neither known whether the
    problem is in~\NP\ nor whether it is in~$\coNP$.\smallskip
    
    The following question is interesting in connection with
    Problem~\ref{prob:facsimp} (see also the notes there): Given an
    (abstract) graph $G$, is $G$ the graph of a simple polytope? If we
    do not restrict the question to simple polytopes the problem is
    also interesting.
  \end{Comment}
  \Related \ref{prob:facsimp}, \ref{prob:steinitz_problem}
\end{Problem}

\mysection{Beyond Polytopes}
\label{sec:beyond}

This section is concerned with problems on finite abstract simplicial
complexes.  Some of the problems listed are direct generalizations of
problems on polytopes.  Most of the basic notions relevant in our
context can be looked up in~\cite{Zie95}; for topological concepts
like \emph{homology} we refer to Munkres' book~\cite{Mun84}.

A \emph{finite abstract simplicial complex}~$\Delta$ is a non-empty
set of subsets (the \emph{simplices} or \emph{faces}) of a finite set
of \emph{vertices} such that $F\in\Delta$ and $G\subset F$ imply
$G\in\Delta$.  The \emph{dimension} of a simplex $F\in\Delta$ is
$\card{F}-1$. The \emph{dimension}~$\dim(\Delta)$ of~$\Delta$ is the
largest dimension of any of the simplices in~$\Delta$. If all its
maximal simplices with respect to inclusion (i.e., its \emph{facets})
have the same cardinality, then~$\Delta$ is \emph{pure}. A pure
$d$-dimensional finite abstract simplicial complex whose \emph{dual
  graph} (defined on the facets, where two facets are adjacent if they
share a common $(d-1)$-face) is connected is a \emph{pseudo-manifold}
if every $(d-1)$-dimensional simplex is contained in at most two
facets.  The boundary of a simplicial $(d+1)$-dimensional polytope
induces a $d$-dimensional pseudo-manifold.

Throughout this section a finite abstract simplicial complex $\Delta$
is given by its list of facets or by the complete list of all
simplices. In the first case, the input size can be measured by~$n$
and~$m$, the numbers of vertices and facets.

\begin{Problem}{Euler Characteristic}
  \label{prob:euler_characteristic}
  \begin{Definition}
    \Input Finite abstract simplicial complex $\Delta$ given by a list
    of facets%
    \Output Euler characteristic $\chi(\Delta) \in \Z$
  \end{Definition}
  \STATUS{Open}{Polynomial time}
  \begin{Comment}
    It is unknown whether the decision version ``$\chi(\Delta) = 0$?''
    of this problem is in \NP. The problem is easy if $\Delta$ is
    given by a list of all of its simplices. For fixed dimension, one
    can enumerate all simplices of $\Delta$ and compute the Euler
    characteristic in polynomial time.
    
    Currently the fastest way to compute the Euler characteristic is
    to determine $\mathcal{V} = \{ S : S$ is an intersection of facets
    of $\Delta\}$ and then compute $\chi(\Delta)$ in time
    $\order{\card{\mathcal{V}}^2}$ by a M\"{o}bius function approach,
    see Rota~\cite{Rot64}. Usually $\mathcal{V}$ is much smaller than
    the whole face lattice of $\Delta$. $\mathcal{V}$ can be listed
    lexicographically by an algorithm of Ganter~\cite{Gan87}, in time
    $\order{\min\{m,n\} \cdot \alpha \cdot \card{\mathcal{V}}}$, where
    $\alpha$ is the number of vertex-facets incidences.
  \end{Comment}
  \Related \ref{prob:sc_f_vector}
\end{Problem}

\begin{ProblemIndex}{$f$-Vector of Simplicial Complexes}{F-vector of simplicial complexes}
  \label{prob:sc_f_vector}
  \begin{Definition}
    \Input Finite abstract simplicial complex $\Delta$ given by a list
    of facets%
    \Output The $f$-vector of $\Delta$
  \end{Definition}
  \STATUS{\NumP-hard}{Polynomial time}
  \begin{Comment}
    If $\Delta$ is given by all of its simplices the problem is
    trivial. Clearly, for fixed $k$, the first $k$ entries of the
    $f$-vector can be computed in polynomial time, since the number of
    $k$-simplices in $\Delta$ is polynomial in $n$. Hence the problem
    is polynomial time solvable for fixed dimension $\dim(\Delta)$.
    
    It is unknown whether the decision problem ``Given the 
    list of facets of~$\Delta$ and some $\varphi\in\N$; is~$\varphi$
    the total number of faces of~$\Delta$?'' is contained in~\NP. This
    problem is only known to be in~\NP\ for partitionable (see
    Problem~\ref{prob:USO}) simplicial complexes (see Kleinschmidt and
    Onn~\cite{KleO96}).\smallskip
    
    To the best of our knowledge, no proof of \NumP-hardness of the
    general problem has appeared in the literature. Therefore we
    include one here.
    
    Consider an instance of SAT, i.e., a formula in conjunctive normal
    form (CNF-formula) $C_1 \wedge \dots \wedge C_m$ with variables
    $x_1, \dots, x_n$ (each $C_i$ contains only disjunctions of
    literals). It is well known (Valiant~\cite{Val79}) that computing
    the number of satisfying truth assignments is \NumP-complete.
    Define $E = \{t_1, f_1, \dots, t_n, f_n\}$.
    
    Part I. First, let $E$ be the vertex set of a simplicial complex
    $\Delta$ defined by the minimal non-faces (\emph{circuits}) $C_1',
    \dots, C_m', P_1, \dots, P_n$, where $P_i = \{t_i, f_i\}$ for
    every $i$.  Here for any clause $C$, $C' := \{ f_j : x_j$ literal
    in $C\} \cup \{ t_j : \overline{x}_j$ literal in $C\}$, e.g., for
    $C = x_1 \vee x_2 \vee \overline{x}_3$ we have $C' = \{f_1, f_2,
    t_3\}$.  The idea is that $t_i$ corresponds to the assignment of a
    true-value and $f_i$ corresponds to the assignment of a
    false-value to variable $x_i$.  The circuits exclude subsets of
    $E$ which include both $t_i$ and $f_i$ for all variables $x_i$ and
    exclude truth-assignments to variables which would not satisfy a
    clause $C_j$. It is, however, allowed that for some variable $x_i$
    neither $t_i$ nor $f_i$ is included in a face. But every
    $(n-1)$-face ($n$-subset of $E$) (if there exists any) corresponds
    to a truth-assignment to the variables (which uses exactly one
    value for each variable) and satisfies the formula.  These subsets
    are counted by $f_{n-1}(\Delta)$. Hence computing $f_{n-1}$ is
    \NumP-complete and computing the $f$-vector of $\Delta$ is
    \NumP-hard. Moreover this shows that computing the dimension of a
    simplicial complex given by the minimal non-faces is
    \NP-hard.
    
    Part II. We now construct a simplicial complex $\overline{\Delta}$
    (the \emph{dual} complex) which is given by facets. Define
    $\overline{\Delta}$ by the facets $\overline{C_1'}, \dots,
    \overline{C_m'}, \overline{P_1}, \dots, \overline{P_n}$, where for
    $S \subseteq E$, $\overline{S} := E \setminus S$. We have that a
    set $S \subseteq E$ is a face of $\Delta$ if and only if
    $\overline{S}$ is \emph{not} a face of $\overline{\Delta}$. Hence,
    $f_{n-1}(\Delta) + f_{n-1}(\overline{\Delta}) = \binom{2n}{n}$.
    It follows that one can efficiently compute $f_{n-1}(\Delta)$ from
    $f_{n-1}(\overline{\Delta})$.
  \end{Comment}
  \Related \ref{prob:f_vector}, \ref{prob:euler_characteristic}
\end{ProblemIndex}

\begin{Problem}{Homology}
  \begin{Definition}
    \Input Finite abstract simplicial complex $\Delta$ given by a list
    of facets, $i\in\N$%
    \Output The $i$-th homology group of $\Delta$, given by its rank
    and its torsion coefficients
  \end{Definition}
  \STATUS{Open}{Polynomial time}
  \begin{Comment}
    There exists a polynomial time algorithm if $\Delta$ is given by
    the list of all simplices, since the Smith normal form of an
    integer matrix can be computed efficiently
    (Iliopoulos~\cite{Ili89}). For fixed $i$ or $\dim(\Delta)-i$, the
    sizes of the boundary matrices are polynomial in the size of
    $\Delta$ and the Smith normal form can again be computed
    efficiently.
  \end{Comment}
  \Related \ref{prob:euler_characteristic}, \ref{prob:sc_f_vector}
\end{Problem}

\begin{Problem}{Shellability}
  \label{prob:shellability}
  \begin{Definition}
    \Input Finite abstract pure simplicial complex $\Delta$ given by a
    list
    of facets%
    \Output ``Yes'' if $\Delta$ is shellable, ``No'' otherwise
  \end{Definition}
  \STATUS{Open}{Open}
  \begin{Comment}
    Given an ordering of the facets of $\Delta$, it can be tested in
    polynomial time whether it is a shelling order. Hence, the problem
    in \NP.
    
    The problem can be solved in polynomial time for one-dimensional
    complexes, i.e., for graphs: a graph is shellable if and only if
    it is connected. Even for $\dim(\Delta)=2$, the status is open. In
    particular, it is unclear if the problem can be solved in
    polynomial time if $\Delta$ is given by a list of all simplices.
    
    For two-dimensional pseudo-manifolds the problem can be solved in
    linear time (Danarj and Klee~\cite{DanK78}).
  \end{Comment}
  \Related \ref{prob:AOF}, \ref{prob:partitionability}
\end{Problem}

\begin{Problem}{Partitionability}
  \label{prob:partitionability}
  \begin{Definition}
    \Input Finite abstract simplicial complex $\Delta$ given by a list
    of facets%
    \Output ``Yes'' if $\Delta$ is partionable, ``No'' otherwise
  \end{Definition}
  \STATUS{Open}{Open}
  \begin{Comment}
    As in Problem \ref{prob:USO}, partitionability is meant in the sense
    of Stanley~\cite{Sta79} (see also~\cite{Zie95}).
    Noble~\cite{Nob96} proved that the problem is in \NP.
    
    \textsc{Partitionability} can be solved in polynomial time for
    one-dimensional complexes, i.e., for graphs: a graph is
    partitionable if and only if at most one of its connected
    components is a tree.  For two-dimensional complexes the
    complexity status is open. In particular, it is unclear if the
    problem can be solved in polynomial time if $\Delta$ is given by a
    list of all simplices.
  \end{Comment}
  \Related \ref{prob:USO}, \ref{prob:shellability}
\end{Problem}

\newpage
\noindent\textbf{\large Table of Problems}\medskip
\begingroup
\setcounter{tocdepth}{2}
\tableofproblems
\endgroup
\newpage

\end{document}